\documentclass[a4paper,12pt,titlepage]{article}
\usepackage{epsfig}
\usepackage{amsmath,amsfonts,amssymb}
\usepackage{graphicx}
\usepackage{epsfig}

\begin{document}

\begin{center}
{\LARGE Quantile regression when the covariates are functions}

\vspace{0.5cm}

Herv\'e Cardot\footnote{INRA Toulouse, Biom\'etrie et Intelligence Artificielle, 31326 Castanet-Tolosan Cedex, France}, Christophe Crambes\footnote{Universit\'e Paul Sabatier, Laboratoire de Statistique et Probabilit\'es, UMR C5583, 118, route de Narbonne, 31062 Toulouse Cedex, France} and Pascal Sarda\footnote{Universit\'e Toulouse-le-Mirail, GRIMM, EA 3686, 5, all\'ees Antonio Machado, 31058 Toulouse Cedex 9, France}

\vspace{0.5cm}

{\footnotesize E-mail addresses: \footnotemark[1]cardot@toulouse.inra.fr, \footnotemark[2]crambes@cict.fr, \footnotemark[3]sarda@cict.fr}

\vspace{0.5cm}

\end{center}

\paragraph{Abstract}

This paper deals with a linear model of regression on quantiles when the explanatory variable takes values in some functional space and the response is scalar. We propose a spline estimator of the functional coefficient that minimizes a penalized $L^{1}$ type criterion. Then, we study the asymptotic behavior of this estimator. The penalization is of primary importance to get existence and convergence. 

\paragraph{Key words} Functional data analysis, conditional quantiles, $B$-spline functions, roughness penalty.

\section{Introduction}

Because of the increasing performances of measurement apparatus and computers, many data are collected and saved on thinner and thinner time scales or spatial grids (temperature curves, spectrometric curves, satellite images, \ldots). So, we are led to process data comparable to curves or more generally to functions of continuous variables (time, space). These data are called \textit{functional data} in the literature (see Ramsay and Silverman, 2002). Thus, there is a need to develop statistical procedures as well as theory for this kind of data and actually many recent works study models taking into account the functional nature of the data. 

Mainly in a formal way, the oldest works in that direction intended to give a mathematical framework based on the theory of linear operators in Hilbert spaces (see Deville, 1974, Dauxois and Pousse, 1976). After that and in an other direction, practical aspects of extensions of descriptive statistical methods like for example Principal Component Analysis have been considered (see Besse and Ramsay, 1986). The monographs by Ramsay and Silverman (1997, 2002) are important contributions in this area.

As pointed out by Ramsay and Silverman (1997), ``{\it the goals of functional data analysis are essentially the same as those of other branches of Statistics}'': one of this goal is the explanation of variations of a dependent variable $Y$ (response) by using information from an independent functional variable $X$ (explanatory variable). In many applications, the response is a scalar: see Frank and Friedman (1993), Ramsay and Silverman (1997), ... Traditionally, one deals, for such a problem, with estimating the regression on the mean  {\it i.e.} the minimizer among some class of functionals $r$ of

$$ \mathbb{E} \left[ (Y-r(X))^2 \right]. $$

As when $X$ is a vector of real numbers, the two main approaches are linear (see Ramsay and Dalzell, 1991, for the functional linear model) or purely nonparametric (see Ferraty and Vieu, 2002, which adapt kernel estimation to the functional setting). It is also known that estimating the regression on the median or more generally on quantiles has some interest. The problem is then to estimate the minimizer among $g_\alpha$ of

\begin{eqnarray}\label{eq1}
\mathbb{E} \left[ l_{\alpha} \left( Y - g_{\alpha} (X) \right) \right],
\end{eqnarray}
where $l_\alpha(u)=|u|+(2\alpha-1)u$. The value $\alpha=1/2$ corresponds to the conditional median whereas values $\alpha\in]0,1[$ correspond to conditional quantiles of order $\alpha$. The advantage of estimating conditional quantiles may be found in many applications such as in agronomy (estimation of yield thresholds), in medicine or in reliability. Besides robust aspects of the median, it may also help to derive some kind of confidence prediction intervals based on quantiles. 

In our work, we assume that the conditional quantile of order $\alpha$ can be written as

\begin{eqnarray}\label{eq2}
g_{\alpha} (X) = \langle \Psi_{\alpha},X \rangle,
\end{eqnarray}
where $<.,.>$ is a functional inner product and the parameter of the model $\Psi_\alpha$ is a function to be estimated. This is the equivalent of the linear model for regression quantiles studied by Koenker and Bassett (1978) where the inner product is the Euclidean one and the parameter is a vector of scalars. We choose to estimate the function $\Psi_\alpha$ by a ``direct'' method: writing our estimator as a linear combination of B-splines, it minimizes the empirical version of expectation (\ref{eq1}) with the addition of a penalty term proportional to the square norm of a given order derivative of the spline. The penalization term allows on one side to control the regularity of the estimator and on the other side to get consistency.

Unlike for the square function, minimization of function $l_\alpha$ does not lead to an explicit expression of the estimator. While computation of the estimator can be resolved by using traditional algorithms (for instance based on Iteratively Weighted Least Squares), the convexity of $l_\alpha$ allows theoretical developments.

In section \ref{sec2}, we define more precisely the framework of our study and the spline estimator of the functional parameter $\Psi_\alpha$. Section \ref{sec3} is devoted to the asymptotic behaviour of our estimator: we study $L^2$ convergence and derive an upper bound for the rate of convergence. Comments on the model and on the optimality of the rate of convergence are given in section \ref{sec4}. Finally, the  proofs are gathered in section \ref{sec5}.


\section{Construction of the estimator}\label{sec2}

In this work, the data consist of an i.i.d. sample of pairs $(X_{i},Y_{i})_{i=1, \ldots, n}$ drawn from a population distribution $(X,Y)$. 
We consider explanatory variables $X_i$ which are square integrable (random) functions defined on $[0,1]$, {\it i.e.} are elements of the space $L^{2} ([0,1])$ so that $X_i=(X_i(t),t\in[0,1])$. The response $Y_i$ is a scalar belonging to $\mathbb{R}$. Assume that $H$, the range of $X$, is a closed subspace of $L^{2} ([0,1])$. For $Y$ having a finite expectation, $\mathbb{E} (|Y|) <+\infty$, and for $\alpha\in]0,1[$, the \textit{conditional $\alpha$-quantile} functional $g_{\alpha}$ of $Y$ given $X$ is a functional defined on $H$ minimizing (\ref{eq1}).

Our aim is to generalize the linear model introduced by Koenker and Bassett (1978). In our setting, it consists in assuming that $g_\alpha$ is a linear and continuous functional defined on $H$ and then it follows that $g_\alpha(X)$ can be written as in (\ref{eq2}). Taking the usual inner product in $L^2([0,1])$, we can write

$$ g_{\alpha}(X) = \langle \Psi_{\alpha},X \rangle = \int_{0}^{1} \Psi_{\alpha} (t) X (t) \; dt, $$ 

\noindent where $\Psi_\alpha$ is the functional coefficient in $H$ to be estimated, the order $\alpha$ being fixed. From now on we consider, for simplicity, that the random variables  $X_{i}$ are centered, that is to say $\mathbb{E} (X_{i}(t)) = 0$, for $t$ a. e.

When $X$ is multivariate, Bassett and Koenker (1978) study the {\it least absolute error (LAE)} estimator for the conditional median, which can be extended to any quantile replacing the absolute value by the convex function $l_\alpha$ in the criterion to be minimized (see Koenker and Bassett, 1978). In our case where we have to estimate a function belonging to an infinite dimensional space, we are looking at an estimator in the form of an expansion in some basis of $B$-splines functions and then minimizing a similar criterion with however the addition of a penalty term.

Before describing in details the estimation procedure, let us note that estimation of conditional quantiles has received a special attention in the multivariate case. As said before, linear modelling has been mainly investigated by Bassett and Koenker (1978). For nonparametric models, we may distinguish  two different approaches: ``indirect'' estimators which are based on a preliminary estimation of the conditional cumulative distribution function (cdf) and ``direct'' estimators which are based on the minimizing the empirical version of criterion (\ref{eq1}). In the class of ``indirect'' estimators, Bhattacharya and Gangopadhyay (1990) study a kernel estimator of the conditional cdf, and estimation of the quantile is achieved by inverting this estimated cdf. In the class of ``direct'' estimators, kernel estimators based on local fit have been proposed (see Tsybakov, 1986, Lejeune and Sarda, 1988 or Fan, Hu and Truong, 1994); in a similar approach, He and Shi (1994) and Koenker \textit{et. al.} (1994) propose a spline estimator. Although our setting is quite different, we adapt in our proofs below some arguments of the work by  He and Shi (1994).

In nonparametric estimation, it is usual to assume that the function to be estimated is sufficiently smooth so that it can be expended in some basis: the degree of smoothness is quantified by the number of derivatives and a lipschitz condition for the derivative of greatest order (see condition (H.2) below). It is also quite usual to approximate such kind of functions by means of regression splines (see de Boor, 1978, for a guide for splines). For this, we have to select a degree $q$ in $\mathbb{N}$ and a subdivision of $[0,1]$ defining the position of the knots. Although it is not necessary, we take equispaced knots so that only the number of the knots has to be selected: for $k$ in $\mathbb{N}^{\star}$, we consider $k-1$ knots that define a subdivision of the interval $[0,1]$ into $k$ sub-intervals. For asymptotic theory, the degree $q$ is fixed but the number of sub-intervals $k$ depends on the sample size $n,$ $k=k_n.$ 
It is well-known that a spline function is a piecewise polynomial: we consider here piecewise polynomials of degree $q$ on each sub-interval, and $(q-1)$ times differentiable on $[0,1]$. This space of spline functions is a vectorial space of dimension $k+q$. 
A basis of this vectorial space is the set of the so-called normalized $B$-spline functions, that we note by $\mathbf{B}_{k,q} = \; \left( B_{1}, \ldots, B_{k+q} \right)^\tau$.

Then, we estimate $\Psi_{\alpha}$ by a linear combination of functions $B_{l}$. This leads us to find a vector $\boldsymbol{\widehat{\theta}} = ( \widehat{\theta}_{1}, \ldots, \widehat{\theta}_{k+q} )^{\tau}$ in $\mathbb{R}^{k+q}$ such that

\begin{equation}\label{estimator} \widehat{\Psi}_{\alpha} = \sum_{l=1}^{k+q} \widehat{\theta}_{l} B_{l} = \mathbf{B}_{k,q}^{\tau} \boldsymbol{\widehat{\theta}}.
\end{equation}
\noindent It is then natural to look for $\widehat{\Psi}_{\alpha}$ as the minimizer of the empirical version of (\ref{eq1}) among functional $g_\alpha$ of the form (\ref{eq2}) with functions $\Psi_\alpha$ belonging to the space of spline functions defined above. 
We will however consider a penalized criterion as we will see now. In our setting, the pseudo-design matrix $\mathbf{A}$ is the matrix of dimension $n\times(k+q)$ and elements $\langle X_{i} , B_{j} \rangle$ for $i = 1, \ldots, n$ and $j=1, \ldots, k+q$. Even if we do not have an explicit expression for a solution to the minimization problem, it is known that the solution would depend on the properties of the inverse of the matrix $\frac{1}{n} \mathbf{A}^{\tau} \mathbf{A}$ which is the $(k+q) \times (k+q)$ matrix with general term $\langle \Gamma_{n} (B_{j}),B_{l} \rangle$, where $\Gamma_{n}$ is the empirical version of the covariance operator $\Gamma_{X}$ of $X$ defined for all $u$ in $L^{2}([0,1])$ by 

\begin{equation}
\Gamma_{X} u = \mathbb{E} \left( \langle X,u \rangle X \right). \label{gammax}
\end{equation}

\noindent We know that $\Gamma_{X}$ is a nuclear operator (see Dauxois \textit{et al}, 1982), consequently no bounded inverse exists for this operator. Moreover, as a consequence of the first monotonicity principle (see theorem 7.1, p.58, in Weinberger, 1974), the restriction of this operator to the space of spline functions has smaller eigenvalues than $\Gamma_{X}$. Finally, it appears to be impossible to control the speed of convergence to zero of the smallest eigenvalue of $\frac{1}{n} \mathbf{A}^{\tau} \mathbf{A}$ (when $n$ tends to infinity): in that sense, we are faced with an inversion problem that can be qualified as ill-conditioned. 
A way to circumvent this problem is to introduce a penalization term in the minimization criterion (see Ramsay and Silverman, 1997, or Cardot \textit{et al.}, 2003, for a similar approach in the functional linear model). Thus, the main role of the penalization is to control the inversion of the matrix linked to the solution of the problem and it consists in restricting the space of solutions. The penalization introduced below will have another effect since we also want to control the smoothness of our estimator. For this reason, and following several authors (see references above), we choose a penalization which allows to control the  norm of the derivative of order $m>0$ of any linear combination  of B-spline functions, so that it can be expressed matricially. Denoting by $(\mathbf{B}_{k,q}^{\tau} \boldsymbol{\theta})^{(m)}$ the $m$-th derivative of the spline function $\mathbf{B}_{k,q}^{\tau} \boldsymbol{\theta},$ we have
$$  \left\| ( \mathbf{B}_{k,q}^{\tau} \boldsymbol{\theta})^{(m)} \right\|^{2} = \boldsymbol{\theta}^{\tau} \mathbf{G}_{k} \boldsymbol{\theta}, \quad \forall \boldsymbol{\theta} \in \mathbb{R}^{k+q},
$$
where $\mathbf{G}_{k}$ is the $(k+q) \times (k+q)$ matrix with general term $[\mathbf{G}_{k}]_{jl} = \langle B_{j}^{(m)},B_{l}^{(m)} \rangle.$

 Then, the vector $\boldsymbol{\widehat{\theta}}$ in (\ref{estimator}) is chosen as the solution of the following minimization problem

\begin{equation}
\min_{\boldsymbol{\theta} \in \mathbb{R}^{k+q}} \biggl\lbrace \frac{1}{n} \sum_{i=1}^{n} l_{\alpha} (Y_{i} -  \langle \mathbf{B}_{k,q}^{\tau} \boldsymbol{\theta},X_{i} \rangle) + \rho \parallel (\mathbf{B}_{k,q}^{\tau} \boldsymbol{\theta})^{(m)} \parallel^{2}
\biggr\rbrace, \label{minfct}
\end{equation}

\noindent where  $\rho$ is the penalization parameter. In the next section, we present a convergence result of the solution of (\ref{minfct}). Note that the role of the penalization also clearly appears in this result.


\section{Convergence result} \label{sec3}

We present in this section the main result on the convergence of our estimator. 
The behaviour of our estimator is linked to a penalized version of the matrix $\mathbf{\widehat{C}}=\frac{1}{n} \mathbf{A}^{\tau} \mathbf{A}$. 
More precisely, adopting the same notations as in Cardot \textit{et. al.} (2003), the existence and convergence of our estimator depend on the inverse of the matrix $\mathbf{\widehat{C}}_{\rho} = \mathbf{\widehat{C}} + \rho \mathbf{G}_{k}$. Under the hypotheses of theorem \ref{thm1} below, the smallest eigenvalue of $\mathbf{\widehat{C}}_{\rho}$, noted $\lambda_{\textrm{min}} (\mathbf{\widehat{C}}_{\rho})$, tends to zero as the sample size $n$ tends to infinity. As the rate of convergence of $\widehat{\Psi}_{\alpha}$ depends on the speed of convergence of $\lambda_{\textrm{min}} (\mathbf{\widehat{C}}_{\rho})$ to zero, we introduce a sequence $(\eta_{n})_{n \in \mathbb{N}}$ such that the set $\Omega_{n}$ defined by

\begin{equation}
\Omega_{n} = \left\{ \omega / \lambda_{\textrm{min}} (\mathbf{\widehat{C}}_{\rho}) > c \eta_{n} \right\}, \label{omega}
\end{equation}

\noindent has probability which goes to 1 when $n$ goes to infinity. Cardot \textit{et al.} (2003) have shown that such a sequence exists in the sense that under hypotheses of theorem \ref{thm1}, there exists a strictly positive sequence $(\eta_{n})_{n \in \mathbb{N}}$ tending to zero as $n$ tends to infinity and such that 

\begin{equation}
\lambda_{\textrm{min}} (\mathbf{\widehat{C}}_{\rho}) \geq c \eta_{n} + o_{P} \left( (k_{n}^{2} n^{1 - \delta})^{-1/2} \right), \label{lambdamin}
\end{equation}

\noindent with $\delta \in ]0,1[$.

To prove the convergence result of the estimator $\widehat{\Psi}_{\alpha}$, we assume that the following hypotheses are satisfied.

$(H.1)$ $\parallel X \parallel \leq C_{0} < + \infty, \quad a.s.$

$(H.2)$ The function $\Psi_{\alpha}$ is supposed to have a $p'$-th derivative $\Psi_{\alpha}^{(p')}$ such that

\vspace{-0.3cm} $$ \quad \left| \Psi_{\alpha}^{(p')} (t) - \Psi_{\alpha}^{(p')} (s) \right| \leq C_{1} | t - s |^{\nu}, \quad s,t \in [0;1],$$

\noindent where $C_{1} > 0$ and $\nu \in [0,1]$. In what follows, we set $p = p' + \nu$ and we suppose that $q \geq p \geq m$.

$(H.3)$ The eigenvalues of $\Gamma_{X}$ (defined in (\ref{gammax})) are strictly positive.

$(H.4)$ For $x\in H$, the random variable $\epsilon$ defined by $\epsilon = Y - \langle \Psi_{\alpha} , X \rangle$ has conditional density function $f_{x}$ given $X=x$, continuous and bounded below by a strictly positive constant at 0, uniformly for $x\in H$.

We derive in theorem \ref{thm1} below an upper bound for the rate of convergence with respect to some kind of $L^2$-norm. Indeed, the operator $\Gamma_X$ is strictly non-negative, so we can associate it a semi-norm noted $\left\| . \right\|_{2}$ and defined by $\left\| u \right\|_{2}^{2} = \langle \Gamma_{X} u,u \rangle$. Then, we have the following result.

\newtheorem{thm}{Theorem}
\newtheorem{cor}{Corollary}

\begin{thm} \label{thm1}
Under hypotheses $(H.1)-(H.4)$, if we also suppose that there exists $\beta, \gamma$ in $]0,1[$ such that $k_{n} \sim n^{\beta}$, $\rho \sim n^{-\gamma}$ and $\eta_{n} \sim n^{-\beta - (1-\delta)/2}$ (where $\delta$ is defined in relation (\ref{lambdamin})), then

(i) $\widehat{\Psi}_{\alpha}$ exists and is unique except on a set whose probability goes to zero as $n$ goes to infinity,

(ii) ${\displaystyle \parallel \widehat{\Psi}_{\alpha} - \Psi_{\alpha} \parallel_{2}^{2}  = O_{P} \left( \frac{1}{k_{n}^{2p}} + \frac{1}{n \eta_{n}} + \frac{\rho^{2}}{k_{n} \eta_{n}} + \rho k_{n}^{2(m-p)} \right).}$
\end{thm}


\section{Some comments} \label{sec4}

(i) Hypotheses (H.1) and (H.3) are quite usual in the functional setting: see for instance Bosq (2000) or Cardot {\it et al.} (2003). Hypothesis $(H.4)$ implies uniqueness of the conditional quantile of order $\alpha$.

\noindent (ii) Some arguments in the proof of theorem \ref{thm1} are inspired from the demonstration of He and Shi (1994) within the framework of real covariates. Moreover, some results from Cardot \textit{et. al.} (2003) are also useful, mainly to deal with the penalization term as pointed out above.\\
Note that it is assumed in the model of He and Shi (1994) that the error term is independent of $X$: condition (H.4) allows us to deal with a more general setting, as in Koenker and Bassett (1978).

\noindent (iii) It is possible to choose particular values for $\beta$ and $\gamma$ to optimize the upper bound for the rate of convergence in theorem \ref{thm1}. In particular, we remark the importance to control the speed of convergence to 0 of the smallest eigenvalue of $\mathbf{\widehat{C}}_{\rho}$ by $\eta_{n}$. For example, Cardot \textit{et al.} (2003) have shown that, under hypotheses of theorem \ref{thm1}, relation (\ref{lambdamin}) is true with $\eta_{n} = \rho / k_{n}$. This gives us

$$ \parallel \widehat{\Psi}_{\alpha} - \Psi_{\alpha} \parallel_{2}^{2} = O_{P} \left( \frac{1}{k_{n}^{2p}} + \frac{k_{n}}{n \rho} + \rho + \rho k_{n}^{2(m-p)} \right). $$

\noindent A corollary is obtained if we take $k_{n} \sim n^{1 / (4p+1)}$ and $\rho \sim n^{-2p / (4p+1)}$; then we get

$$\parallel \widehat{\Psi}_{\alpha} - \Psi_{\alpha} \parallel_{2}^{2} = O_{P} \left(  n^{-2p / (4p+1)} \right). $$

\noindent We can imagine that, with stronger hypotheses on the random function $X$, we can find a sequence $\eta_{n}$ greater than $\rho / k_{n}$, that will improve the convergence speed of the estimator. As a matter of fact, the rate derived in theorem \ref{thm1} does not imply the rate obtained by Stone (1982), that is to say a rate of order $n^{-2p/(2p+1)}$. Indeed, suppose that $1 / k_{n}^{2p}$, $1 / (n \eta_{n})$ and $\rho^{2} / (k_{n} \eta_{n})$ are all of order $n^{-2p/(2p+1)}$. This would imply that $k_{n} \sim n^{1 / (2p+1)}$ and $\eta_{n} \sim n^{-1 / (2p+1)}$, which contradicts the condition $\eta_{n} \sim n^{-\beta - (1-\delta)/2}$. Nevertheless, it is possible to obtain a speed of order $n^{-2p/(2p+1) + \kappa}$. This leads to $k_{n} \sim n^{1 / (2p+1) - \kappa / (2p)}$ and $\eta_{n} \sim n^{-1 / (2p+1) - \kappa}$. Then, the condition $\eta_{n} \sim n^{-\beta - (1-\delta)/2}$ implies $\kappa = p(1 - \delta) / (2p+1)$. So finally, we get $k_{n} \sim n^{(1 + \delta) / 2(2p+1)}$, $\rho \sim n^{(-4p -1 + \delta) / 4(2p+1)}$ and $\eta_{n} \sim n^{(-p -1 +p \delta) / (2p+1)}$. The convergence result would be then

$$ \parallel \widehat{\Psi}_{\alpha} - \Psi_{\alpha} \parallel_{2}^{2} = O_{P} \left(  n^{-p(1 + \delta) / (2p+1)} \right). $$

\noindent A final remark is that the last term $\rho k_{n}^{2(m-p)}$ of the speed in theorem \ref{thm1} is not always negligible compared to the other terms. However, it will be the case if we suppose that $m \leq p / (1 + \delta) + (1 - \delta) / 4 (1 + \delta)$.

\noindent (iv) This quantile estimator is quite useful in practice, specially for forecasting purpose (by conditional median or inter-quantiles intervals). From a computational point of view, several algorithms may be used: we have implemented in the $R$ language an algorithm based on the Iterated Reweighted Least Square (IRLS). Note that even for real data cases, the curves are always observed in some discretization points, the regression splines is easy to implement by approximating inner products with quadrature rules. The IRLS algorithm (see Ruppert and Carroll, 1988, Lejeune and Sarda, 1988) allows to build conditional quantiles spline estimators and gives satisfactory forecast results. This algorithm has been used in particular on  the ``ORAMIP'' (``Observatoire R\'egional de l'Air en Midi-Pyr\'en\'ees'') data to forecast pollution in the city of Toulouse (France): the results of this practical study are described in Cardot \textit{et. al.} (2004). We are interested in predicting the ozone concentration one day ahead, 
knowing the ozone curve (concentration along time) the day before. In that special case, conditional quantiles were also useful to predict an ozone threshold such that the probability to exceed this threshold is a given risk $1-\alpha$. In other words, it comes back to give an estimation of the $\alpha$-quantile maximum ozone knowing the ozone curve the day before.


\section{Proof of theorem \ref{thm1}} \label{sec5}

The proof of the result is based on the same kind of decomposition of $\widehat{\Psi}_{\alpha} - \Psi_{\alpha}$ as the one used by He and Shi (1994). The main difference comes from the fact that our design matrix is ill-conditioned, which led us to add the penalization term treated using some arguments from Cardot \textit{et al.} (2003).

\vspace{0.5cm}

\noindent Hypothesis $(H.2)$ implies (see de Boor, 1978) that there exists a spline function $\Psi_{\alpha}^{\star} = \mathbf{B}_{k,q}^{\tau} \boldsymbol{\theta^{\star}}$, called \textit{spline approximation} of $\Psi_{\alpha}$, such that

\begin{equation}
\sup_{t \in [0,1]} \left| \Psi_{\alpha}^{\star} (t) - \Psi_{\alpha} (t) \right| \leq \frac{C_{2}}{k_{n}^{p}}. \label{approx}
\end{equation}

\noindent In what follows, we set $R_{i} = \langle \Psi_{\alpha}^{\star} - \Psi_{\alpha} , X_{i} \rangle$; so we deduce from (\ref{approx}) and from hypothesis $(H.1)$ that there exists a positive constant $C_{3}$ such that

\begin{equation}
\max_{i=1, \ldots, n} | R_{i} | \leq \frac{C_{3}}{k_{n}^{p}}, \quad a.s.
\label{reste}
\end{equation}

\noindent The operator $\Gamma_{n}$ allows to define the empirical version of the $L^{2}$ norm by $\left\| u \right\|_{n}^{2} = \langle \Gamma_{n} u , u \rangle$. At first, we show the result (ii) of theorem \ref{thm1} for the penalized empirical $L^{2}$ norm. Writing $\widehat{\Psi}_{\alpha} - \Psi_{\alpha} = (\widehat{\Psi}_{\alpha} - \Psi_{\alpha}^{\star}) + (\Psi_{\alpha}^{\star} - \Psi_{\alpha})$, we get

\begin{eqnarray*}
& & \| \widehat{\Psi}_{\alpha} - \Psi_{\alpha} \|_{n}^{2} + \rho \| (\widehat{\Psi}_{\alpha} - \Psi_{\alpha})^{(m)} \|^{2} \\ 
& \leq & \frac{2}{n} \sum_{i=1}^{n} \langle \widehat{\Psi}_{\alpha} - \Psi_{\alpha}^{\star} , X_{i} \rangle^{2} + \frac{2}{n} \sum_{i=1}^{n} \langle \Psi_{\alpha}^{\star} - \Psi_{\alpha} , X_{i} \rangle^{2} \\ 
& & + 2 \rho \| (\widehat{\Psi}_{\alpha} - \Psi_{\alpha}^{\star})^{(m)} \|^{2} + 2 \rho \| (\Psi_{\alpha}^{\star} - \Psi_{\alpha})^{(m)} \|^{2}.
\end{eqnarray*}

\noindent Now, using again hypothesis $(H.1)$, we get almost surely and for all $i = 1, \ldots, n$, the inequality $\langle \Psi_{\alpha}^{\star} - \Psi_{\alpha} , X_{i} \rangle^{2} \leq C_{0}^{2} C_{2}^{2} / k_{n}^{2p}$. Moreover, lemma 8 of Stone (1985) gives us the existence of a positive constant $C_{4}$ that satisfies $\| (\Psi_{\alpha} - \Psi_{\alpha}^{\star})^{(m)} \|^{2} \leq C_{4} k_{n}^{2(m-p)}$. So we deduce

\begin{eqnarray}
& & \| \widehat{\Psi}_{\alpha} - \Psi_{\alpha} \|_{n}^{2} + \rho \| (\widehat{\Psi}_{\alpha} - \Psi_{\alpha})^{(m)} \|^{2} \nonumber \\
& \leq & \frac{2}{n} \sum_{i=1}^{n} \langle \widehat{\Psi}_{\alpha} - \Psi_{\alpha}^{\star} , X_{i} \rangle^{2} + 2 \rho \| (\widehat{\Psi}_{\alpha} - \Psi_{\alpha}^{\star})^{(m)} \|^{2}  \nonumber \\
& & + \frac{2 C_{0}^{2} C_{2}^{2}}{k_{n}^{2p}} + 2 C_{4} \rho k_{n}^{2(m-p)}, \quad a.s. \label{empirique}
\end{eqnarray}

\noindent Our goal is now to compare our estimator $\widehat{\Psi}_{\alpha}$ with the spline approximation $\Psi_{\alpha}^{\star}$. For that, we adopt the following transformation $\boldsymbol{\theta} = \mathbf{\widehat{C}}_{\rho}^{-1/2} \boldsymbol{\beta} + \boldsymbol{\theta^{\star}}$. Then, we define on the set $\Omega_{n}$

\begin{eqnarray*}
f_{i} (\boldsymbol{\beta}) & = & l_{\alpha} \left[ Y_{i} - \langle \mathbf{B}_{k,q}^{\tau} \left( \mathbf{\widehat{C}}_{\rho}^{-1/2} \boldsymbol{\beta} + \boldsymbol{\theta^{\star}} \right) , X_{i} \rangle \right] \\
& & + \rho \left\| \left[ \mathbf{B}_{k,q}^{\tau} \left( \mathbf{\widehat{C}}_{\rho}^{-1/2} \boldsymbol{\beta} + \boldsymbol{\theta^{\star}} \right) \right]^{(m)} \right\|^{2}.
\end{eqnarray*}

\noindent We notice that minimizing $\sum_{i=1}^{n} f_{i} (\boldsymbol{\beta})$ comes back to the minimization of the criterion (\ref{minfct}). We are interested by the behaviour of the function $f_{i}$ around zero: $f_{i} (\boldsymbol{0})$ is the value of our loss criterion when $\boldsymbol{\theta} = \boldsymbol{\theta^{\star}}$. Let us also notice that the inverse of the matrix $\mathbf{\widehat{C}}_{\rho}$ appears in the definition of $f_{i}$. This inverse exists on the set $\Omega_{n}$ defined by (\ref{omega}), and which probability goes to $1$ as $n$ goes to infinity. Lemma \ref{lem1} below, whose proof is given in section \ref{prooflem1}, allows us to get the results (i) and (ii) of theorem \ref{thm1} for the penalized empirical $L^{2}$ norm.

\newtheorem{lemme}{Lemma}

\begin{lemme} \label{lem1}
Under the hypotheses of theorem \ref{thm1}, for all $\epsilon > 0$, there exists $L = L_{\epsilon}$ (sufficiently large) and $(\delta_{n})_{n \in \mathbb{N}}$ with $\delta_{n} = \sqrt{1/(n \eta_{n}) + \rho^{2}/(k_{n} \eta_{n})}$ such that, for $n$ large enough

$$ P \left[ \inf_{|\boldsymbol{\beta}|=L \delta_{n}} \sum_{i=1}^{n} f_{i} (\boldsymbol{\beta}) > \sum_{i=1}^{n} f_{i} (\mathbf{0}) \right] > 1 - \epsilon. $$
\end{lemme}

\noindent Using convexity arguments, this inequality means that the solution $\boldsymbol{\widehat{\beta}}$ exists and is unique on the ball centered in $\boldsymbol{\theta^{\star}}$ and of radius $L \delta_{n}$. As we use the one-to-one transformation $\boldsymbol{\theta} = \mathbf{\widehat{C}}_{\rho}^{-1/2} \boldsymbol{\beta} + \boldsymbol{\theta^{\star}}$ on the set $\Omega_{n}$, we deduce the existence and the uniqueness of the solution of (\ref{minfct}) on the set $\Omega_{n}$, which proves point (i) of theorem \ref{thm1}.

\vspace{0.2cm}

\noindent Now, let $\epsilon$ be strictly positive; using the convexity of function $f_{i}$, there exists $L = L_{\epsilon}$ such that, for $n$ large enough

\begin{equation}
P \left[ \inf_{|\boldsymbol{\beta}| \geq L \delta_{n}} \sum_{i=1}^{n} f_{i} (\boldsymbol{\beta}) > \sum_{i=1}^{n} f_{i} (\mathbf{0}) \right] > 1 - \epsilon. \label{res1}
\end{equation}

\noindent On the other hand, using the definition of $f_{i}$ and the minimization criterion (\ref{minfct}), we have

\begin{eqnarray*}
& & \frac{1}{n} \sum_{i=1}^{n} f_{i} \left( \mathbf{\widehat{C}}_{\rho}^{1/2} \boldsymbol{\widehat{\theta}} - \mathbf{\widehat{C}}_{\rho}^{1/2} \boldsymbol{\theta^{\star}} \right) \\
& = & \inf_{\boldsymbol{\theta} \in \mathbb{R}^{k+q}} \left[ \frac{1}{n} \sum_{i=1}^{n} l_{\alpha} \left( Y_{i} - \langle \mathbf{B}_{k,q}^{\tau} \boldsymbol{\theta} , X_{i} \rangle \right) + \rho \left\| \left( \mathbf{B}_{k,q}^{\tau} \boldsymbol{\theta} \right)^{(m)} \right\|^{2} \right], 
\end{eqnarray*}

\noindent so we finally get

$$ \frac{1}{n} \sum_{i=1}^{n} f_{i} \left( \mathbf{\widehat{C}}_{\rho}^{1/2} \widehat{\boldsymbol{\theta}} - \mathbf{\widehat{C}}_{\rho}^{1/2} \boldsymbol{\theta^{\star}} \right) \leq \frac{1}{n} \sum_{i=1}^{n} f_{i} (\mathbf{0}). $$

\noindent Then, combining this with equation (\ref{res1}), we obtain

\begin{equation}
P \left[ \inf_{|\boldsymbol{\beta}| \geq L \delta_{n}} \sum_{i=1}^{n} f_{i} (\boldsymbol{\beta}) > \sum_{i=1}^{n} f_{i} \left( \mathbf{\widehat{C}}_{\rho}^{1/2} \boldsymbol{\widehat{\theta}} - \mathbf{\widehat{C}}_{\rho}^{1/2} \boldsymbol{\theta^{\star}} \right) \right] > 1 - \epsilon. \label{res2}
\end{equation}

\noindent Now, using the definition of $\mathbf{\widehat{C}}_{\rho}$, we have

\begin{eqnarray*}
& & P \Bigg[ \frac{1}{n} \sum_{i=1}^{n} \langle \widehat{\Psi}_{\alpha} - \Psi_{\alpha}^{\star} , X_{i} \rangle^{2} + \rho \left\| (\widehat{\Psi}_{\alpha} - \Psi_{\alpha}^{\star})^{(m)} \right\|^{2} \leq L^{2} \delta_{n}^{2} \Bigg] \\
& = & 1 - P \left[ \left| \mathbf{\widehat{C}}_{\rho}^{1/2} (\boldsymbol{\widehat{\theta}} - \boldsymbol{\theta^{\star}}) \right| > L \delta_{n} \right] \\
& \geq & P \left[ \inf_{|\boldsymbol{\beta}| \geq L \delta_{n}} \sum_{i=1}^{n} f_{i} (\boldsymbol{\beta}) > \sum_{i=1}^{n} f_{i} \left( \mathbf{\widehat{C}}_{\rho}^{1/2} \boldsymbol{\widehat{\theta}} - \mathbf{\widehat{C}}_{\rho}^{1/2} \boldsymbol{\theta^{\star}} \right) \right].
\end{eqnarray*}

\noindent With relation (\ref{res2}), this last probability is greater than $1 - \epsilon$, so we obtain

$$ \frac{1}{n} \sum_{i=1}^{n} \langle \widehat{\Psi}_{\alpha} - \Psi_{\alpha}^{\star} , X_{i} \rangle^{2} + \rho \left\| (\widehat{\Psi}_{\alpha} - \Psi_{\alpha}^{\star})^{(m)} \right\|^{2} = O_{P} \left( \delta_{n}^{2} \right) = O_{P} \left( \frac{1}{n \eta_{n}} + \frac{\rho^{2}}{k_{n} \eta_{n}} \right). $$

\vspace{0.2cm}

\noindent This last result, combined with inequality (\ref{empirique}) finally gives us the equivalent of result (ii) for the penalized empirical $L^{2}$ norm. Point (ii) ( with the norm $\left\| . \right\|_{2}$) then follows from lemma \ref{lem2} below, which is proved in section \ref{prooflem2}, and achieves the proof of theorem \ref{thm1} (ii).

\begin{lemme} \label{lem2}
Let $f$ and $g$ be two functions supposed to be $m$ times differentiable and such that

\vspace{-0.3cm}

$$ \| f - g \|_{n}^{2} + \rho \| \left( f - g \right)^{(m)} \|^{2} = O_{P} ( u_{n} ), $$

\noindent with $u_{n}$ going to zero when $n$ goes to infinity. Under hypotheses $(H.1)$ and $(H.3)$ and if moreover $\| g \|$ and $\| g^{(m)} \|$ are supposed to be bounded, we have

\vspace{-0.3cm}

$$ \| f - g \|_{2}^{2} = O_{P} ( u_{n} ). $$
\end{lemme}

\subsection{Proof of lemma \ref{lem1}} \label{prooflem1}

This proof is based on three preliminary lemmas, proved respectively in sections \ref{prooflem3}, \ref{prooflem4} and \ref{prooflem5}. We denote by $T_{n}$ the set of the random variables $(X_{1}, \ldots, X_{n})$. Under hypotheses of theorem \ref{thm1}, we have the following results.

\begin{lemme} \label{lem3}
There exists a constant $C_{5}$ such that, on the set $\Omega_{n}$ defined by (\ref{omega}), we have

$$ \max_{i=1, \ldots, n} \left| \langle \mathbf{B}_{k,q}^{\tau} \mathbf{\widehat{C}}_{\rho}^{-1/2} \boldsymbol{\beta} , X_{i} \rangle \right| \leq \frac{C_{5} | \boldsymbol{\beta} |}{\sqrt{k_{n} \eta_{n}}}, \quad a.s. $$
\end{lemme}

\begin{lemme} \label{lem4}
For all $\epsilon > 0$, there exists $L = L_{\epsilon}$ such that

$$ \lim_{n \to +\infty} P \left[ \inf_{| \boldsymbol{\beta} | = 1} \sum_{i=1}^{n} \left( f_{i} (L \delta_{n} \boldsymbol{\beta}) - f_{i} (\mathbf{0}) - \mathbb{E} \left[ f_{i} (L \delta_{n} \boldsymbol{\beta}) - f_{i} (\mathbf{0}) | T_{n} \right] \right) > \epsilon \delta_{n}^{2} n \right] = 0. $$
\end{lemme}

\begin{lemme} \label{lem5}
For all $\epsilon > 0$, there exists $L = L_{\epsilon}$ such that

$$ P \left[ \inf_{| \boldsymbol{\beta} | = 1} \sum_{i=1}^{n} \mathbb{E} \left[ f_{i} (L \delta_{n} \boldsymbol{\beta}) - f_{i} (\mathbf{0}) | T_{n} \right] > \delta_{n}^{2} n \right] > 1 - \epsilon. $$
\end{lemme}

\noindent These three lemmas allow us to prove lemma \ref{lem1}. Indeed, let $L$ be a strictly positive real number; we make the following decomposition 

$$ \inf_{|\boldsymbol{\beta}|=1} \sum_{i=1}^{n} f_{i} (L \delta_{n} \boldsymbol{\beta}) - \sum_{i=1}^{n} f_{i} (\mathbf{0}) \geq A_{n} + B_{n}, $$

\noindent with 

$$ A_{n} = \inf_{|\boldsymbol{\beta}|=1} \sum_{i=1}^{n} \left( f_{i} (L \delta_{n} \boldsymbol{\beta}) - f_{i} (\mathbf{0}) - \mathbb{E} \left[ f_{i} (L \delta_{n} \boldsymbol{\beta}) - f_{i} (\mathbf{0}) | T_{n} \right] \right) $$ 

\noindent and 

$$ B_{n} = \inf_{|\boldsymbol{\beta}|=1} \sum_{i=1}^{n} \mathbb{E} \left[ f_{i} (L \delta_{n} \boldsymbol{\beta}) - f_{i} (\mathbf{0}) | T_{n} \right]. $$

\noindent Using lemmas \ref{lem4} and \ref{lem5}, we can find $L$ sufficiently large such that, for $n$ large enough

$$ P \left( \left| A_{n} \right| > \delta_{n}^{2} n  \right) < \epsilon / 2, $$

\noindent and

$$ P \left( B_{n} > \delta_{n}^{2} n \right) > 1 - \epsilon / 2, $$  

\noindent thus we get

$$ P \left[ \inf_{|\boldsymbol{\beta}|=1} \sum_{i=1}^{n} f_{i} (L \delta_{n} \boldsymbol{\beta}) - \sum_{i=1}^{n} f_{i} (\mathbf{0}) > 0 \right] \geq P \left( A_{n} + B_{n} > 0 \right) > 1 - \epsilon, $$

\noindent which achieves the proof of lemma \ref{lem1}.

\subsection{Proof of lemma \ref{lem3}} \label{prooflem3}

Using lemma 6.2 of Cardot \textit{et al.} (2003), we have 

$$ \lambda_{\textrm{min}} (\mathbf{\widehat{C}}_{\boldsymbol{\rho}}) \geq C_{5}' \eta_{n} + o_{P} ( ( k_{n}^{2} n^{1-\delta} )^{-1/2} ). $$ 

\noindent Noticing that $\left| \langle \mathbf{B}_{k,q}^{\tau} \mathbf{\widehat{C}}_{\rho}^{-1/2} \boldsymbol{\beta} , X_{i} \rangle \right|^{2} \leq \langle \mathbf{B}_{k,q}^{\tau} , X_{i} \rangle \mathbf{\widehat{C}}_{\rho}^{-1} \langle \mathbf{B}_{k,q} , X_{i} \rangle | \boldsymbol{\beta} |^{2}$, we deduce that

\begin{eqnarray*}
& & \left| \langle \mathbf{B}_{k,q}^{\tau} \mathbf{\widehat{C}}_{\rho}^{-1/2} \boldsymbol{\beta} , X_{i} \rangle \right|^{2} \\
& \leq & \langle \mathbf{B_{k,q}^\tau} , X_{i} \rangle \langle \mathbf{B_{k,q}} , X_{i} \rangle | \boldsymbol{\beta} |^{2} \left[ \frac{1}{C_{5}' \eta_{n}} + o_{P} \left( \left( k_{n}^{2} n^{1-\delta} \right)^{-1/2} \right) \right], 
\end{eqnarray*}

\noindent which gives us $\left| \langle \mathbf{B_{k,q}^\tau} \mathbf{\widehat{C}}_{\boldsymbol{\rho}}^{-1/2} \boldsymbol{\beta} , X_{i} \rangle \right|^{2} \leq C_{5}'' | \boldsymbol{\beta} |^{2} / (k_{n} \eta_{n}) + o_{P} \left( n^{(\delta-1)/2} \right)$ almost surely, and achieves the proof of lemma \ref{lem3}.

\subsection{Proof of lemma \ref{lem4}} \label{prooflem4}

Considering the definition of functions $f_{i}$ and $l_{\alpha}$, we have 

\begin{eqnarray*}
& & \sup_{| \boldsymbol{\beta} | \leq 1} \sum_{i=1}^{n} \Big( f_{i} (L \delta_{n} \boldsymbol{\beta}) - f_{i} (\mathbf{0}) - \mathbb{E} \left[ f_{i} (L \delta_{n} \boldsymbol{\beta}) - f_{i} (\mathbf{0}) | T_{n} \right] \Big) \\ 
& = & \sup_{| \boldsymbol{\beta} | \leq 1} \sum_{i=1}^{n} \Big( \Big| \epsilon_{i} - L \delta_{n} \langle \;  \mathbf{B}_{k,q}^{\tau} \mathbf{\widehat{C}}_{\rho}^{-1/2} \boldsymbol{\beta} , X_{i} \rangle - R_{i} \Big| - | \epsilon_{i} - R_{i} | \\
& & - \mathbb{E} \Big[ \Big| \epsilon_{i} - L \delta_{n} \langle \mathbf{B}_{k,q}^{\tau} \mathbf{\widehat{C}}_{\rho}^{-1/2} \boldsymbol{\beta} , X_{i} \rangle - R_{i} \Big| - | \epsilon_{i} - R_{i} | | T_{n} \Big] \Big),
\end{eqnarray*}

\noindent where $\epsilon_{1}, \ldots, \epsilon_{n}$ are $n$ real random variables independent and identically distributed defined by $\epsilon_{i} = Y_{i} - \langle \Psi_{\alpha} , X_{i} \rangle$ for all $i=1, \ldots, n$. Let us also denote $\Delta_{i} (\boldsymbol{\beta}) = \Big| \epsilon_{i} - L \delta_{n} \langle \mathbf{B}_{k,q}^{\tau} \mathbf{\widehat{C}}_{\rho}^{-1/2} \boldsymbol{\beta} , X_{i} \rangle - R_{i} \Big| - | \epsilon_{i} - R_{i} |$. To prove lemma \ref{lem4}, it suffices to show that, for all $\epsilon > 0$, there exists $L = L_{\epsilon}$ such that

$$ \lim_{n \to +\infty} P \left( \sup_{| \boldsymbol{\beta} | \leq 1} \sum_{i=1}^{n} \left[ \Delta_{i} (\boldsymbol{\beta}) - \mathbb{E} ( \Delta_{i} (\boldsymbol{\beta}) | T_{n}) \right] > \epsilon \delta_{n}^{2} n \right) = 0. $$

\noindent Let $\epsilon$ be a real number strictly positive and $\mathcal{C}$ the subset of $\mathbb{R}^{k+q}$ defined by $\mathcal{C} = \left\{ \boldsymbol{\beta} \in \mathbb{R}^{k+q} / | \boldsymbol{\beta} | \leq 1 \right\}$. As $\mathcal{C}$ is a compact set, we can cover it with open balls, that is to say $\mathcal{C} = \bigcup_{j=1}^{K_{n}} \mathcal{C}_{j}$ with $K_{n}$ chosen, for all $j$ from $1$ to $K_{n}$, such that

\begin{equation}
\textrm{diam } (\mathcal{C}_{j}) \leq \frac{\epsilon \delta_{n} \sqrt{k_{n} \eta_{n}}}{8 C_{5} L}. \label{cj}
\end{equation} 

\noindent Hence

\begin{equation}
K_{n} \leq \left( \frac{8 C_{5} L}{\epsilon \delta_{n} \sqrt{k_{n} \eta_{n}}} \right)^{k_{n}+q}. \label{compact}
\end{equation} 

\noindent Now, for $1 \leq j \leq K_{n}$, let $\boldsymbol{\beta_{j}}$ be in $\mathcal{C}_{j}$; using the definition of $\Delta_{i} (\boldsymbol{\beta})$ and the triangular inequality, we have

\begin{eqnarray*}
& & \min_{j=1, \ldots, K_{n}} \sum_{i=1}^{n} \left| \left[ \Delta_{i} (\boldsymbol{\beta}) - \mathbb{E} ( \Delta_{i} (\boldsymbol{\beta}) | T_{n}) \right] - \left[ \Delta_{i} (\boldsymbol{\beta_{j}}) - \mathbb{E} ( \Delta_{i} (\boldsymbol{\beta}_{j}) | T_{n}) \right] \right| \\
& \leq & 2 L \delta_{n} \min_{j=1, \ldots, K_{n}} \sum_{i=1}^{n} \left| \langle \mathbf{B}_{k,q}^{\tau} \mathbf{\widehat{C}}_{\rho}^{-1/2} (\boldsymbol{\beta} - \boldsymbol{\beta_{j}}) , X_{i} \rangle \right|.
\end{eqnarray*}

\noindent Then, using lemma \ref{lem3}, we get

\begin{eqnarray*}
& & \min_{j=1, \ldots, K_{n}} \sum_{i=1}^{n} \left| \left[ \Delta_{i} (\boldsymbol{\beta}) - \mathbb{E} ( \Delta_{i} (\boldsymbol{\beta}) | T_{n}) \right] - \left[ \Delta_{i} (\boldsymbol{\beta_{j}}) - \mathbb{E} ( \Delta_{i} (\boldsymbol{\beta}_{j}) | T_{n}) \right] \right| \\
& \leq & 2 L \delta_{n} \frac{C_{5} n}{\sqrt{k_{n} \eta_{n}}} \min_{j=1, \ldots, K_{n}} | \boldsymbol{\beta} - \boldsymbol{\beta_{j}} |, 
\end{eqnarray*}

\noindent this last inequality being true only on the set $\Omega_{n}$ defined by (\ref{omega}). Moreover, there exists a unique $j_{0} \in \{ 1, \ldots, K_{n} \}$ such that $\boldsymbol{\beta} \in \mathcal{C}_{j_{0}}$, which gives us with relation (\ref{cj})

\begin{equation}
\min_{j=1, \ldots, K_{n}} \sum_{i=1}^{n} \left| \left[ \Delta_{i} (\boldsymbol{\beta}) - \mathbb{E} ( \Delta_{i} (\boldsymbol{\beta}) | T_{n}) \right] - \left[ \Delta_{i} (\boldsymbol{\beta_{j}}) - \mathbb{E} ( \Delta_{i} (\boldsymbol{\beta}_{j}) | T_{n}) \right] \right| \leq \frac{\epsilon}{4} \delta_{n}^{2} n. \label{ineg1}
\end{equation}

\noindent On the other hand, we have

$$ \sup_{\boldsymbol{\beta} \in \mathcal{C}} | \Delta_{i} (\boldsymbol{\beta}) | \leq L \delta_{n} \sup_{\boldsymbol{\beta} \in \mathcal{C}} | \langle \mathbf{B}_{k,q}^{\tau} \mathbf{\widehat{C}}_{\rho}^{-1/2} \boldsymbol{\beta} , X_{i} \rangle |, $$

\noindent and using lemma \ref{lem3} again, we get, on $\Omega_{n}$,

\begin{equation}
\sup_{\boldsymbol{\beta} \in \mathcal{C}} | \Delta_{i} (\boldsymbol{\beta}) | \leq \frac{C_{5} L \delta_{n}}{\sqrt{k_{n} \eta_{n}}}. \label{Berns1}
\end{equation}

\noindent Besides, for $\boldsymbol{\beta}$ fixed in $\mathcal{C}$, with the same arguments as before, if we denote by $T^{\star}$ the set of the random variables $(X_{1}, \ldots, X_{n}, \ldots)$, we have

$$ \sum_{i=1}^{n} \textrm{Var} \left( \Delta_{i} (\boldsymbol{\beta}) | T^{\star} \right) \leq \sum_{i=1}^{n} L^{2} \delta_{n}^{2} \textrm{Var} \left( | \langle \mathbf{B}_{k,q}^{\tau} \mathbf{\widehat{C}}_{\rho}^{-1/2} \boldsymbol{\beta} , X_{i} \rangle |^{2} | T^{\star} \right). $$

\noindent Then, using the definition of $\mathbf{\widehat{C}}_{\rho}$, we remark that

\begin{equation}
\sum_{i=1}^{n} \Big| \langle \mathbf{B}_{k,q}^{\tau} \mathbf{\widehat{C}}_{\rho}^{-1/2} \boldsymbol{\beta} , X_{i} \rangle \Big|^{2} = n | \boldsymbol{\beta} |^{2} - n \rho  \boldsymbol{\beta}^{\tau} \mathbf{\widehat{C}}_{\rho}^{-1/2} \mathbf{G}_{k} \mathbf{\widehat{C}}_{\rho}^{-1/2} \boldsymbol{\beta}, \label{rel*}
\end{equation}

\noindent which gives us

\begin{equation}
\sum_{i=1}^{n} \textrm{Var} \left( \Delta_{i} (\boldsymbol{\beta}) | T^{\star} \right) \leq n L^{2} \delta_{n}^{2}. \label{Berns2}
\end{equation}

\noindent We are now able to prove lemma \ref{lem4}. Using first relation (\ref{ineg1}), we have

\begin{eqnarray*}
& & P \left[ \left( \sup_{| \boldsymbol{\beta} | \leq 1} \sum_{i=1}^{n} \left[ \Delta_{i} (\boldsymbol{\beta}) - \mathbb{E} \left( \Delta_{i} (\boldsymbol{\beta}) | T_{n} \right)  \right] > \epsilon \delta_{n}^{2} n \right) \cap \Omega_{n} \Big| T^{\star} \right] \\
& \leq & P \left[ \left( \max_{j = 1, \ldots, K_{n}} \sum_{i=1}^{n} \left[ \Delta_{i} (\boldsymbol{\beta_{j}}) - \mathbb{E} \left( \Delta_{i} (\boldsymbol{\beta}_{j}) | T_{n} \right) \right] > \frac{\epsilon}{2} \delta_{n}^{2} n \right) \cap \Omega_{n} \Big| T^{\star} \right],
\end{eqnarray*}

\noindent and then

\begin{eqnarray*}
& & P \left[ \left( \sup_{| \boldsymbol{\beta} | \leq 1} \sum_{i=1}^{n} \left[ \Delta_{i} (\boldsymbol{\beta}) - \mathbb{E} \left( \Delta_{i} (\boldsymbol{\beta}) | T_{n} \right) \right] > \epsilon \delta_{n}^{2} n \right) \cap \Omega_{n} \Big| T^{\star} \right] \\
& \leq & K_{n} P \left[ \left( \sum_{i=1}^{n} \left[ \Delta_{i} (\boldsymbol{\beta_{j}}) - \mathbb{E} \left( \Delta_{i} (\boldsymbol{\beta}_{j}) | T_{n} \right) \right] > \frac{\epsilon}{2} \delta_{n}^{2} n \right) \cap \Omega_{n} \Big| T^{\star} \right].
\end{eqnarray*}

\noindent By inequalities (\ref{Berns1}) and (\ref{Berns2}), we apply Bernstein inequality (see Uspensky, 1937) and inequality (\ref{compact}) to obtain

\begin{eqnarray*}
& & P \left[ \left( \sup_{| \boldsymbol{\beta} | \leq 1} \sum_{i=1}^{n} \left[ \Delta_{i} (\boldsymbol{\beta}) - \mathbb{E} \left( \Delta_{i} (\boldsymbol{\beta}) | T_{n} \right) \right] > \epsilon \delta_{n}^{2} n \right) \cap \Omega_{n} \Big| T^{\star} \right] \\
& \leq & 2 \exp \left\{ \textrm{ ln} \left( \frac{8 C_{5} L n}{\epsilon \delta_{n} \sqrt{k_{n} \eta_{n}}} \right)^{k_{n} + q} - \frac{\epsilon^{2} \delta_{n}^{4} n^{2} / 4}{2 n L^{2} \delta_{n}^{2} + 2 C_{5} L \delta_{n} \times \epsilon \delta_{n}^{2} n / (2 \sqrt{k_{n} \eta_{n}})} \right\}.
\end{eqnarray*}

\noindent This bound does not depend on the sample $T^{\star} = (X_{1}, \ldots, X_{n}, \ldots)$, hence, if we take the expectation on both sides of this inequality above, we deduce

\begin{eqnarray*}
& & P \left[ \left( \sup_{\boldsymbol{\beta} \leq 1} \sum_{i=1}^{n} \left[ \Delta_{i} (\boldsymbol{\beta}) - \mathbb{E} \left( \Delta_{i} (\boldsymbol{\beta}) | T_{n} \right) \right] > \epsilon \delta_{n}^{2} n \right) \cap \Omega_{n} \right] \\
& \leq & 2 \exp \left\{ - \frac{\epsilon^{2} \delta_{n}^{2} \sqrt{k_{n} \eta_{n}} n}{8 L^{2} \sqrt{k_{n} \eta_{n}} + 4 C_{5} L \delta_{n}} \right. \\
& & \left. \times \left[ 1 - \frac{(k_{n} + q) (8 L^{2} \sqrt{k_{n} \eta_{n}} + 4 C_{5} L \delta_{n})}{\epsilon^{2} \delta_{n}^{2} \sqrt{k_{n} \eta_{n}} n} \textrm{ ln} \left( \frac{8 C_{5} L n}{\epsilon \delta_{n} \sqrt{k_{n} \eta_{n}}} \right) \right] \right\}.
\end{eqnarray*}

\noindent If we fix $L = L_{n} = \sqrt{n k_{n} \eta_{n} \delta_{n}^{2}}$, we have

\begin{itemize}
\item[] ${\displaystyle \frac{\delta_{n}^{2} \sqrt{k_{n} \eta_{n}} n}{L^{2} \sqrt{k_{n} \eta_{n}}} = \frac{1}{k_{n} \eta_{n}} \xrightarrow[n \to +\infty]{} +\infty}$,
\item[]  
\item[] ${\displaystyle \frac{\delta_{n}^{2} \sqrt{k_{n} \eta_{n}} n}{L \delta_{n}} = \sqrt{n} \xrightarrow[n \to +\infty]{} +\infty}$, 
\item[] 
\item[] ${\displaystyle \frac{k_{n} L^{2} \sqrt{k_{n} \eta_{n}}}{\delta_{n}^{2} \sqrt{k_{n} \eta_{n}} n} = k_{n}^{2} \eta_{n} \xrightarrow[n \to +\infty]{} 0}$,
\item[] 
\item[] ${\displaystyle \frac{k_{n} L \delta_{n}}{\delta_{n}^{2} \sqrt{k_{n} \eta_{n}} n} = \frac{k_{n}}{\sqrt{n}} \xrightarrow[n \to +\infty]{} 0}$. 
\item[] 
\end{itemize}

\noindent This leads to

$$ \lim_{n \to +\infty} P \left[ \left( \sup_{| \boldsymbol{\beta} | \leq 1} \sum_{i=1}^{n} \left[ \Delta_{i} (\boldsymbol{\beta}) - \mathbb{E} \left( \Delta_{i} (\boldsymbol{\beta}) | T_{n} \right) \right] > \epsilon \delta_{n}^{2} n \right) \cap \Omega_{n} \right] = 0, $$

\noindent and with the fact that $\Omega_{n}$ has probability tending to $1$ when $n$ goes to infinity, we finally obtain

$$ \lim_{n \to +\infty} P \left[ \sup_{| \boldsymbol{\beta} | \leq 1} \sum_{i=1}^{n} \left[ \Delta_{i} (\boldsymbol{\beta}) - \mathbb{E} \left( \Delta_{i} (\boldsymbol{\beta}) | T_{n} \right) \right] > \epsilon \delta_{n}^{2} n \right] = 0, $$

\noindent which achieves the proof of lemma \ref{lem4}.

\subsection{Proof of lemma \ref{lem5}} \label{prooflem5}

Let $a$ and $b$ be two real numbers. We denote by $F_{i \epsilon}$ the random repartition function of $\epsilon_{i}$ given $T_{n}$ and by $f_{i \epsilon}$ the random density function of $\epsilon_{i}$ given $T_{n}$. As $\mathbb{E} \left( l_{\alpha} (\epsilon_{i} + b) | T_{n} \right) = \int_{\mathbb{R}} l_{\alpha} (s + b) \; d F_{i \epsilon} (s)$, we obtain, using a Taylor linearization at first order, the existence of a quantity $r_{iab}$ such that

$$ \mathbb{E} \left( l_{\alpha} (\epsilon_{i}+a+b) - l_{\alpha} (\epsilon_{i}+b) | T_{n} \right) = f_{i \epsilon} (0) a^{2} + 2 f_{i \epsilon} (0) ab + (\frac{a^{2}}{2} + ab) r_{iab}, $$

\noindent with $r_{iab} \longrightarrow 0$ when $a,b \longrightarrow 0$. If we set $L' = \sqrt{2} L$ and $R_{i}' = \sqrt{2} R_{i}$, this relation gives us

\begin{eqnarray}
& & \sum_{i=1}^{n} \mathbb{E} \left[ l_{\alpha} \left( \epsilon_{i} - L \delta_{n} \langle  \mathbf{B}_{k,q}^{\tau} \mathbf{\widehat{C}}_{\rho}^{-1/2} \boldsymbol{\beta}  ,X_{i} \rangle - R_{i} \right) - l_{\alpha} \left( \epsilon_{i} - R_{i} \right) | T_{n} \right] \nonumber \\ 
& = & 2 \sum_{i=1}^{n} f_{i \epsilon} (0) \left[ L'^{2} \delta_{n}^{2} \langle \mathbf{B}_{k,q}^{\tau} \mathbf{\widehat{C}}_{\rho}^{-1/2} \boldsymbol{\beta} ,X_{i} \rangle^{2} + L' \delta_{n} \langle  \mathbf{B}_{k,q}^{\tau} \mathbf{\widehat{C}}_{\rho}^{-1/2} \boldsymbol{\beta}  ,X_{i} \rangle R_{i}' \right] \nonumber \\ 
& & \! \! \! \! \! \! + \sum_{i=1}^{n} \left[ L'^{2} \delta_{n}^{2} \langle  \mathbf{B}_{k,q}^{\tau} \mathbf{\widehat{C}}_{\rho}^{-1/2} \boldsymbol{\beta}  ,X_{i} \rangle^{2} + L' \delta_{n} \langle \mathbf{B}_{k,q}^{\tau} \mathbf{\widehat{C}}_{\rho}^{-1/2} \boldsymbol{\beta}  ,X_{i} \rangle R_{i}' \right] r_{i \boldsymbol{\beta}}, \label{eg1}
\end{eqnarray}

\noindent with $r_{i \boldsymbol{\beta}} \longrightarrow 0$. Considering $\boldsymbol{\beta}$ such that $| \boldsymbol{\beta} | = 1$, we have, using relation (\ref{reste})

\begin{eqnarray}
& & L'^{2} \delta_{n}^{2} \langle  \mathbf{B}_{k,q}^{\tau} \mathbf{\widehat{C}}_{\rho}^{-1/2} \boldsymbol{\beta}  ,X_{i} \rangle^{2} + L' \delta_{n} \langle \mathbf{B}_{k,q}^{\tau} \mathbf{\widehat{C}}_{\rho}^{-1/2} \boldsymbol{\beta} , X_{i} \rangle R_{i}' \nonumber \\ 
& \geq & \frac{1}{2} \; L'^{2} \delta_{n}^{2} \langle \mathbf{B}_{k,q}^{\tau} \mathbf{\widehat{C}}_{\rho}^{-1/2} \boldsymbol{\beta}  ,X_{i} \rangle^{2} - \frac{C_{3}^{2}}{k_{n}^{2p}}, \quad a.s. \label{ineg2}
\end{eqnarray}

\noindent Moreover, if we set $V_{n} = \sup_{| \boldsymbol{\beta} | = 1} \max_{i = 1 \ldots n} | r_{i \boldsymbol{\beta}} |$, then with condition (H.4) $1 \! \! 1_{\{ V_{n} < \min_{i} f_{i \epsilon} (0) / 4 \}} = 1 \! \! 1_{\mathbb{R}}$ for $n$ large enough, and

\begin{eqnarray}
& & \left| \left[ L'^{2} \delta_{n}^{2} \langle \mathbf{B}_{k,q}^{\tau} \mathbf{\widehat{C}}_{\rho}^{-1/2} \boldsymbol{\beta}  ,X_{i} \rangle^{2} + L' \delta_{n} \langle \mathbf{B}_{k,q}^{\tau} \mathbf{\widehat{C}}_{\rho}^{-1/2} \boldsymbol{\beta} , X_{i} \rangle R_{i}' \right] r_{i \boldsymbol{\beta}} \right| \nonumber \\ 
& \leq & \frac{\min_{i} f_{i \epsilon} (0)}{4} \left| L'^{2} \delta_{n}^{2} \langle  \mathbf{B}_{k,q}^{\tau} \mathbf{\widehat{C}}_{\rho}^{-1/2} \boldsymbol{\beta}  ,X_{i} \rangle^{2} + L' \delta_{n} \langle \mathbf{B}_{k,q}^{\tau} \mathbf{\widehat{C}}_{\rho}^{-1/2} \boldsymbol{\beta} , X_{i} \rangle R_{i}' \right| \nonumber \\ 
& \leq & 2 \min_{i} f_{i \epsilon} (0) \left[ \frac{3}{16} \; L'^{2} \delta_{n}^{2} \langle \mathbf{B}_{k,q}^{\tau} \mathbf{\widehat{C}}_{\rho}^{-1/2} \boldsymbol{\beta} ,X_{i} \rangle^{2} + \frac{C_{3}^{2}}{8 k_{n}^{2p}} \right].
\label{ineg3}
\end{eqnarray}

\noindent Using inequalities (\ref{ineg2}) and (\ref{ineg3}), relation (\ref{eg1}) becomes then

\begin{eqnarray*}
& & \sum_{i=1}^{n} \mathbb{E} \left[ l_{\alpha} \left( \epsilon_{i} - L \delta_{n} \langle \mathbf{B}_{k,q}^{\tau} \mathbf{\widehat{C}}_{\rho}^{-1/2} \boldsymbol{\beta}  ,X_{i} \rangle - R_{i} \right) - l_{\alpha} \left( \epsilon_{i} - R_{i} \right) | T_{n} \right] \nonumber \\ 
& \geq & 2 \min_{i} f_{i \epsilon} (0) \left[ \frac{5}{16} \; L'^{2} \delta_{n}^{2} \sum_{i=1}^{n} \langle \mathbf{B}_{k,q}^{\tau} \mathbf{\widehat{C}}_{\rho}^{-1/2} \boldsymbol{\beta}  ,X_{i} \rangle^{2} - \frac{9}{8} \; \frac{C_{3}^{2} n}{k_{n}^{2p}} \right]. 
\end{eqnarray*}

\noindent Now, we come back to the definition of function $f_{i}$ to obtain

\begin{eqnarray*}
& & \frac{1}{\delta_{n}^{2} n} \inf_{| \boldsymbol{\beta} | = 1} \sum_{i=1}^{n} \mathbb{E} \left[ f_{i} ( L \delta_{n} \boldsymbol{\beta}) - f_{i}(\mathbf{0}) | T_{n} \right] \\
& \geq & 2 \min_{i} f_{i \epsilon} (0) \left[ \frac{5 L'^{2}}{16 n} \sum_{i=1}^{n} \langle \mathbf{B}_{k,q}^{\tau} \mathbf{\widehat{C}}_{\rho}^{-1/2} \boldsymbol{\beta}  ,X_{i} \rangle^{2} - \frac{9 C_{3}^{2}}{8 k_{n}^{2p} \delta_{n}^{2}} \right] \\
& & \! \! \! \! \! \! + \rho L^{2} \left\| \left( \mathbf{B}_{k,q}^{\tau} \mathbf{\widehat{C}}_{\rho}^{-1/2} \boldsymbol{\beta} \right)^{(m)} \right\|^{2} + 2 \frac{L \rho}{\delta_{n}} \langle \left( \mathbf{B}_{k,q}^{\tau} \mathbf{\widehat{C}}_{\rho}^{-1/2} \boldsymbol{\beta} \right)^{(m)} , \left( \mathbf{B}_{k,q}^{\tau} \boldsymbol{\theta^{\star}} \right)^{(m)} \rangle.
\end{eqnarray*}

\noindent Reminding that $L'{^2} = 2 L^{2}$ and taking $\xi = \min ( \frac{5}{4} \min_{i} f_{i \epsilon} (0) , 1 )$, we have $\xi > 0$ by hypothesis $(H.4)$ and then

\begin{eqnarray*}
& & \frac{1}{\delta_{n}^{2} n} \inf_{| \boldsymbol{\beta} | = 1} \sum_{i=1}^{n} \mathbb{E} \left[ f_{i} ( L \delta_{n} \boldsymbol{\beta}) - f_{i}(\mathbf{0}) | T_{n} \right] \\
& \geq & \xi L^{2} \inf_{| \boldsymbol{\beta} | = 1} \left[ \frac{1}{n} \sum_{i=1}^{n} \langle \mathbf{B}_{k,q}^{\tau} \mathbf{\widehat{C}}_{\rho}^{-1/2} \boldsymbol{\beta}  ,X_{i} \rangle^{2} + \rho \left\| \left( \mathbf{B}_{k,q}^{\tau} \mathbf{\widehat{C}}_{\rho}^{-1/2} \boldsymbol{\beta} \right)^{(m)} \right\|^{2} \right]\\
& & - \frac{9}{4} \min_{i} f_{i \epsilon} (0) \frac{C_{3}^{2}}{k_{n}^{2p} \delta_{n}^{2}} + \frac{2 L \rho}{\delta_{n}} \langle \left( \mathbf{B}_{k,q}^{\tau} \mathbf{\widehat{C}}_{\rho}^{-1/2} \boldsymbol{\beta} \right)^{(m)} , \left( \mathbf{B}_{k,q}^{\tau} \boldsymbol{\theta^{\star}} \right)^{(m)} \rangle.
\end{eqnarray*}

\noindent Using relation (\ref{rel*}), we get

\begin{eqnarray*}
& & \frac{1}{\delta_{n}^{2} n} \inf_{| \boldsymbol{\beta} | = 1} \sum_{i=1}^{n} \mathbb{E} \left[ f_{i} ( L \delta_{n} \boldsymbol{\beta}) - f_{i}(\mathbf{0}) | T_{n} \right] \\
& \geq & \xi L^{2} - \frac{9}{4} \min_{i} f_{i \epsilon} (0) \frac{C_{3}^{2}}{k_{n}^{2p} \delta_{n}^{2}} + \frac{2 L \rho}{\delta_{n}} \langle \left( \mathbf{B}_{k,q}^{\tau} \mathbf{\widehat{C}}_{\rho}^{-1/2} \boldsymbol{\beta} \right)^{(m)} , \left( \mathbf{B}_{k,q}^{\tau} \boldsymbol{\theta^{\star}} \right)^{(m)} \rangle.
\end{eqnarray*}

\noindent Moreover, for $| \boldsymbol{\beta} | = 1$, the infimum of $\langle \big( \mathbf{B}_{k,q}^{\tau} \mathbf{\widehat{C}}_{\rho}^{-1/2} \boldsymbol{\beta} \big)^{(m)} , \big( \mathbf{B}_{k,q}^{\tau} \boldsymbol{\theta^{\star}} \big)^{(m)} \rangle$ is obtained for $\boldsymbol{\beta} = - \mathbf{\widehat{C}}_{\rho}^{1/2} \boldsymbol{\theta^{\star}} / | \mathbf{\widehat{C}}_{\rho}^{1/2} \boldsymbol{\theta^{\star}} |$. Using the fact that the spline approximation has a bounded $m$-th derivative, we deduce the existence of a constant $C_{9}>0$ such that

$$ \inf_{\mid \boldsymbol{\beta} \mid = 1} \langle \big( \mathbf{B}_{k,q}^{\tau} \mathbf{\widehat{C}}_{\rho}^{-1/2} \boldsymbol{\beta} \big)^{(m)} , \big( \mathbf{B}_{k,q}^{\tau} \boldsymbol{\theta^{\star}} \big)^{(m)} \rangle \geq - \frac{C_{9}}{\sqrt{\eta_{n}}}, $$

\noindent hence we obtain

\begin{eqnarray*}
& & \frac{1}{\delta_{n}^{2} n} \inf_{| \boldsymbol{\beta} | = 1} \sum_{i=1}^{n} \mathbb{E} \left[ f_{i} ( L \delta_{n} \boldsymbol{\beta}) - f_{i}(\mathbf{0}) | T_{n} \right] \\
& \geq & \xi L^{2} - \frac{9}{4} \min_{i} f_{i \epsilon} (0) \frac{C_{3}^{2}}{k_{n}^{2p} \delta_{n}^{2}} - 2 C_{9} \frac{L \rho}{\delta_{n} \sqrt{\eta_{n}}},
\end{eqnarray*}

\noindent that is to say

\begin{eqnarray*}
& & \frac{1}{\delta_{n}^{2} n} \inf_{| \boldsymbol{\beta} | = 1} \sum_{i=1}^{n} \mathbb{E} \left[ f_{i} ( L \delta_{n} \boldsymbol{\beta}) - f_{i}(\mathbf{0}) | T_{n} \right] \\
& \geq & \xi L^{2} \left( 1 - \frac{9 \min_{i} f_{i \epsilon} (0) C_{3}^{2}}{4 \xi L^{2} k_{n}^{2p} \delta_{n}^{2}} - \frac{2 C_{9} \rho}{\xi L \delta_{n} \sqrt{\eta_{n}}} \right).
\end{eqnarray*}

\noindent Reminding that we have fixed $L = L_{n} = \sqrt{n k_{n} \eta_{n} \delta_{n}^{2}}$, we get

\begin{itemize}
\item[] for ${\displaystyle \delta_{n}^{2} \sim \frac{1}{n \eta_{n}}, \textrm{ we have } \frac{1}{L^{2} k_{n}^{2p} \delta_{n}^{2}} \sim \frac{k_{n} \eta_{n}}{n \rho^{4} k_{n}^{2p}} \xrightarrow[n \to +\infty]{} 0}$,
\item[] 
\item[] for ${\displaystyle \delta_{n}^{2} \sim \frac{\rho}{k_{n} \eta_{n}}, \textrm{ we have } \frac{\rho}{L \delta_{n} \sqrt{\eta_{n}}} \sim \frac{\rho \sqrt{n}}{\sqrt{k_{n}}} \xrightarrow[n \to +\infty]{} 0}$.
\end{itemize}

\noindent This leads to

$$ \lim_{n \to +\infty} P \left( \frac{1}{\delta_{n}^{2} n} \inf_{| \boldsymbol{\beta} | = 1} \sum_{i=1}^{n} \mathbb{E} \left[ f_{i} ( L \delta_{n} \boldsymbol{\beta}) - f_{i}(\mathbf{0}) | T_{n} \right] > 1 \right) = 0, $$

\noindent which achieves the proof of lemma \ref{lem5}.

\subsection{Proof of lemma \ref{lem2}} \label{prooflem2}

Writing $\Gamma_{X} = (\Gamma_{X} - \Gamma_{n}) + \Gamma_{n}$, we make the following decomposition

\begin{equation}
\| f - g \|_{2}^{2} = 2 \| \Gamma_{X} - \Gamma_{n} \| \left( \| f \|^{2} + \| g \|^{2} \right) + \| f - g \|_{n}^{2}. \label{et1}
\end{equation}

\noindent Now, let us decompose $f$ as follows $f = P + R$ with $P(t) = \sum_{l=0}^{m-1} \frac{t^{l}}{l !} f^{(l)} (0)$ and $R(t) = \int_{0}^{t} \frac{(t-u)^{m-1}}{(m-1) !} f^{(m)} (u) \; du$. $P$ belongs to the space $\mathcal{P}_{m-1}$ of polynomials of degree at most $m-1$, whose dimension is finite and equal to $m$. Using hypothesis $(H.3)$, there exists a constant $C_{6}>0$ such that
we have $\| P \|^{2} \leq C_{6} \| P \|_{n}^{2}$. Then, we can deduce

\begin{eqnarray}
\| f \|^{2} & \leq & 2 \| P \|^{2} + 2 \| R \|^{2} \nonumber \\ 
& \leq & 2 C_{6} \| P \|_{n}^{2} + 2 \| R \|^{2} \nonumber \\ 
& \leq & 4 C_{6} \| f \|_{n}^{2} + 4 C_{6} \| \Gamma_{n} \| \; \| R \|^{2} + 2 \| R \|^{2}. \label{et2}
\end{eqnarray}

\noindent As $\Gamma_{n}$ is a bounded operator (by hypothesis $(H.1)$), there exists a constant $C_{7}>0$ such that we have $\| \Gamma_{n} \| \leq C_{7}$. Moreover, under Cauchy-Schwarz inequality, there exists a constant $C_{8}>0$ such that $\| R \|^{2} \leq C_{8} \| f^{(m)} \|^{2}$. Relation (\ref{et2}) gives $\| f \|^{2} \leq 4 C_{6} \| f \|_{n}^{2} + \left( 4 C_{6} C_{7} + 2 \right) C_{8} \| f^{(m)} \|^{2}$. Then, if we write $f = (f-g)+g$, we finally deduce

\begin{eqnarray}
\| f \|^{2} & \leq & 8 C_{6} \| f - g \|_{n}^{2} + \left( 8 C_{6} C_{7} + 4 \right) C_{8} \| (f-g)^{(m)} \|^{2} \nonumber \\
& & + 8 C_{6} \| \Gamma_{n} \| \; \| g \|^{2} + \left( 8 C_{6} C_{7} + 4 \right) C_{8} \| g^{(m)} \|^{2}. \label{et3}
\end{eqnarray}

\noindent We have supposed that $\| g \|$ and $\| g^{(m)} \|$ are bounded, so

$$ 8 C_{6} \| \Gamma_{n} \| \; \| g \|^{2} + \left( 8 C_{6} C_{7} + 4 \right) C_{8} \| g^{(m)} \|^{2} = O(1), $$

\noindent and the hypothesis $\| f - g \|_{n}^{2} + \rho \| (f-g)^{(m)} \|^{2} = O_{P} (u_{n})$ gives us the bounds $\| f - g \|_{n}^{2} = O_{P} (u_{n})$ and $\| (f-g)^{(m)} \|^{2} = O_{P} \left( u_{n} / \rho \right)$. Then, relation (\ref{et3}) becomes 

\begin{equation}
\| f \|^{2} = O_{P} \Big( 1 + \frac{u_{n}}{\rho} \Big). \label{et4}
\end{equation}

\noindent Finally, we have $\| \Gamma_{X} - \Gamma_{n} \| = o_{P} (n^{(\delta - 1)/2}) = o_{P} (\rho)$ from lemma 5.3 of Cardot \textit{et al.} (1999). This equality, combined with equations (\ref{et1}) and (\ref{et4}) gives us $\| f - g \|_{2}^{2} = O_{P} (u_{n})$, which is the announced result.

\vspace*{0.3cm}

\noindent{\bf Acknowledgements}~: We would like to thank the referees for a careful reading of a previous version of the manuscript that permits to improve some results as well as  the participants to the working group ``staph'' in Toulouse for stimulating discussions.

\section*{Bibliography}

\begin{description}

\item Bassett, G.W. and Koenker, R.W. (1978).
\newblock Asymptotic Theory of Least Absolute Error Regression. 
\newblock {\it J. Amer. Statist. Assoc.}, \textbf{73}, 618-622.

\item Besse, P.C. and Ramsay, J.O. (1986).
\newblock Principal Component Analysis of Sampled Curves.
\newblock \textit{Psychometrika}, \textbf{51}, 285-311.

\item Bhattacharya, P.K. and Gangopadhyay, A.K. (1990).
\newblock Kernel and Nearest-Neighbor Estimation of a Conditional Quantile.
\newblock \textit{Ann. Statist.}, \textbf{18}, 1400-1415.

\item Bosq, D. (2000).
\newblock {\it Linear Processes in Function Spaces}.
\newblock Lecture Notes in Statistics, {\bf 149}, Springer-Verlag.

\item Cardot, H., Ferraty, F. and Sarda, P. (1999).
\newblock The Functional Linear Model.
\newblock \textit{Stat. and Prob. Letters}, \textbf{45}, 11-22.

\item Cardot, H., Ferraty, F. and Sarda, P. (2003).
\newblock Spline Estimators for the Functional Linear Model.
\newblock \textit{Statistica Sinica}, \textbf{13}, 571-591.

\item Cardot, H., Crambes, C. and Sarda, P. (2004).
\newblock Conditional Quantiles with Functional Covariates: an Application to Ozone Pollution Forecasting.
\newblock In \textit{Compstat 2004 Proceedings}, ed. J. Antoch, Physica-Verlag, 769-776.

\item Dauxois, J. and Pousse, A. (1976).
\newblock Les Analyses Factorielles en Calcul des Probabilit\'es et en Statistique: Essai d'\'etude synth\'etique (in French).
\newblock Th\`ese, Universit\'e Paul sabatier, Toulouse, France.

\item Dauxois, J., Pousse, A. and Romain, Y. (1982).
\newblock Asymptotic Theory for the Principal Component Analysis of a Random Vector Function: some Applications to Statistical Inference.
\newblock \textit{J. of Mult. Analysis}, \textbf{12}, 136-154.  

\item Deville, J.C. (1974).
\newblock M\'ethodes Statistiques et Num\'eriques de l'Analyse Harmonique (in French).
\newblock \textit{Ann. Insee}, \textbf{15}.

\item de Boor, C. (1978).
\newblock \textit{A Practical Guide to Splines}.
\newblock Springer, New-York.

\item Fan, J., Hu, T.C. and Truong, Y.K. (1994).
\newblock Robust Nonparametric Function Estimation.
\newblock \textit{Scand. J. Statist}, \textbf{21}, 433-446.

\item Ferraty, F. and Vieu, P. (2002).
\newblock The Functional Nonparametric Model and Application to Spectrometric Data.
\newblock \textit{Computational Statistics}, \textbf{17}, 545-564.

\item He, X. and Shi, P. (1994).
\newblock Convergence Rate of $B$-Spline Estimators of Nonparametric Conditional Quantile Functions.
\newblock \textit{Nonparametric Statistics}, \textbf{3}, 299-308.

\item Koenker, R.W. and Bassett G.W.(1978).
\newblock Regression Quantiles.
\newblock \textit{Econometrica}, \textbf{46}, 33-50.

\item Koenker, R.W., Ng, P. and Portnoy, S. (1994).
\newblock Quantile Smoothing Splines.
\newblock \textit{Biometrika}, \textbf{81}, 673-680.

\item Lejeune, M. and Sarda, P. (1988).
\newblock Quantile Regression: A Nonparametric Approach.
\newblock \textit{Computational Statistics and Data Analysis}, \textbf{6}, 229-239.

\item Ramsay, J.O. and Silverman, B.W. (1997).
\newblock  \textit{Functional Data Analysis}.
\newblock Springer-Verlag.

\item Ramsay, J.O. and Silverman, B.W. (2002).
\newblock  \textit{Applied Functional Data Analysis}.
\newblock Springer-Verlag.

\item Ruppert, D. and Carroll, R.J. (1988).
\newblock \textit{Transformation and Weighting in Regression}.
\newblock Chapman and Hall.

\item Stone, C. (1982).
\newblock Optimal Rates of Convergence for Nonparametric Models.
\newblock \textit{Ann. Statist.}, \textbf{10}, 1040-1053.

\item Stone, C. (1985).
\newblock Additive Regression and other Nonparametric Models.
\newblock \textit{Ann. Statist.}, \textbf{13}, 689-705.

\item Tsybakov, A.B. (1986).
\newblock Robust Reconstruction of Functions by the Local-Approximation Method.
\newblock \textit{Problems of Information Transmission}, \textbf{22}, 133-146.

\item Uspensky, J.V. (1937).
\newblock \textit{Introduction to Mathematical Probability}.
\newblock New York and London.

\item Weinberger, H.F. (1974).
\newblock \textit{Variational Methods for Eigenvalue Approximation}.
\newblock SIAM, Philadelphia.

\end{description}

\end{document}